\documentclass[preprint,12pt]{elsarticle}

\usepackage{amssymb}
\usepackage{amsthm}
\usepackage{amsmath}
\usepackage{graphics}
\usepackage{epsfig}
\usepackage{amssymb}
\usepackage{hyperref}
\usepackage{graphicx}
\usepackage{subfigure}
\usepackage{color}
\usepackage{tikz}
\usetikzlibrary{patterns}
\usepackage[ruled,vlined]{algorithm2e}
\usepackage{hyperref}
\usepackage{a4wide}
\usepackage{amsmath}
\usepackage{amsfonts}
\usepackage{enumerate}
\newcommand{\pf}{\noindent {\bf Proof: }}

\newtheorem*{theoremaux}{Theorem \theoremauxnum}
\gdef\theoremauxnum{1}

\newtheorem{lemma}{\bf Lemma}[section]

\newtheorem{theorem}{\bf Theorem}[section]
\newtheorem{proposition}[lemma]{\bf Proposition}
\newtheorem{corollary}[lemma]{\bf Corollary}
\newtheorem{definition}{\bf Definition}[section]
\newtheorem{remark}{\bf Remark}[section]
\newtheorem{problem}{\bf Problem}
\journal{~}

\begin{document}
	
	\begin{frontmatter}
		
		
		%%\author[rvt]{C.V. ̃Radhakrishnan\corref{cor1}\fnref{fn1}}
%\author{Tapa Manna\dag , Angsuman Das\corref{cor1}, Baby Bhattacharya\dag }
%		\ead{mannatapa24@gmail.com}
%		\address{\dag Department of Mathematics,\\ NIT Agartala, India} \\
 %           \address{Department of Mathematics,\\ Presidency University, Kolkata, India} 
%		\cortext[cor1]{Corresponding author}
 % 
 %		\ead{angsuman.maths@presiuniv.ac.in}
 %		\ead{babybhatt75@gmail.com}

\author[add1]{Tapa Manna}
  \ead{mannatapa24@gmail.com}
  \author[add2]{Angsuman Das\corref{cor1}}
  \ead{angsuman.maths@presiuniv.ac.in}
  \author[add1]{Baby Bhattacharya}
  \ead{babybhatt75@gmail.com}

  \cortext[cor1]{Correspondence Author}
  \address[add1]{Department of Mathematics,\\ NIT Agartala, India}
  \address[add2]{Department of Mathematics,\\ Presidency University, Kolkata, India}

\title{Forbidden Subgraphs of Prime Order Element Graph}

\begin{abstract}
	In this paper, we study different forbidden subgraph characterizations of the prime-order element graph $\Gamma(G)$ defined on a finite group $G$. Its set of vertices is the group $G$ and two vertices $x,y \in G$ are adjacent if the order of $xy$ is prime. More specifically, we investigate the conditions when $\Gamma(G)$ is perfect, cograph, chordal, claw-free, and interval graph.
\end{abstract}
		
\end{frontmatter}

\section{Introduction}
The main motivation for defining graphs on groups (see \cite{Cameron-survey} for an extensive survey) is to study the ambient group from its graph. Forbidden subgraphs have been recently used to characterize groups from their underlying graphs, e.g., \cite{brachter-kaja-1}, \cite{brachter-kaja-2}, \cite{pallabi-cograph}, \cite{pallabi-forbidden}. In this paper, we study some forbidden subgraph characterizations of one such graph, namely the prime-order element graph $\Gamma(G)$ of a group $G$, introduced in \cite{tapa-1}.  We investigate when $\Gamma(G)$ is perfect, cograph, chordal, claw-free, and interval graph. For terms and definitions related to graph theory, one can refer to \cite{west-graph-book}.

\begin{definition} Let $G$ be a finite group. The prime-order element graph $\Gamma(G)$ of a group $G$ is defined to be a graph with $G$ as the set of vertices and two distinct vertices $x$ and $y$ are adjacent if and only if order of $xy$ is prime.	
\end{definition}

We first recall some basic properties of $\Gamma(G)$ from \cite{tapa-1} which will be used in this work.

\begin{proposition}\label{previous-results}
    Let $G$ be a finite group, and let $S$ be the set of all elements of prime order in $G$. Then
    \begin{itemize}
        \item Any vertex of $\Gamma(G)$ have degree either $|S|$ or $|S|-1$.
        \item $\Gamma(G)$ is connected if and only if $\langle S\rangle=G$.
        \item $\Gamma(G)$ has a universal vertex if and only if all nonidentity elements of $G$ are of prime order.  
        \item If $H$ is a subgroup of $G$. Then $\Gamma(H)$ is an induced subgraph of $\Gamma(G)$.
    \end{itemize}
\end{proposition}

Next, we recall a family of finite groups, known as {\it EPPO groups}, which will be used to characterize cographs and clawfree graphs among $\Gamma(G)$. 

\begin{definition} \cite{eppo-shi-yang}, \cite{eppo-suzuki}
    A finite group $G$ is called an EPPO group (or CP group in some old papers) if every nonidentity element of $G$ is of prime power order.
\end{definition}
These are exactly those groups for which the Gruenberg-Kegel graph is edgeless. Equivalently, these are precisely the groups for which the power graph and the enhanced power graph coincide. For a nice classification of EPPO groups, see Theorem 1.7 \cite{eppo-cameron}.

Finally, we state a celebrated theorem, namely the Strong Perfect Graph Theorem, which will be extensively used in this article.

\begin{theorem}[Strong Perfect Graph Theorem \cite{strong-perfect}]
A graph is perfect if and only if it contains no odd hole and no odd antihole of length $\geq 5$.   
\end{theorem}
\subsection{Organization of the paper}
In Section \ref{perfect}, we study the perfectness of $\Gamma(G)$. We start with some well-known families of groups and check whether their corresponding graphs are perfect. Then we provide some necessary and sufficient conditions for a group to have their corresponding graphs perfect. We also study their perfectness based on the order of the underlying groups. In Section \ref{cograph}, we discuss when $\Gamma(G)$ is a cograph. In Section \ref{chordal-section}, we treat chordal graphs, interval graphs, split graphs, and threshold graphs among $\Gamma(G)$. Section \ref{clawfree-section} deals with clawfreeness of $\Gamma(G)$. Finally, we conclude with some open questions in Section \ref{conclusion}, as a future scope in this direction.

\section{When is $\Gamma(G)$ perfect}\label{perfect}
In this section, we study the perfectness of $\Gamma(G)$. To start with, we focus on some well-known family of groups regarding whether they yield perfect graphs.
\begin{theorem}\label{bunch}
The following are true:
\begin{enumerate}
    \item \label{pq-not-perfect} $\Gamma(\mathbb{Z}_{pq})$ where $p,q\neq 2$ is not perfect.
    \item \label{2p*2-not-perfect} $\Gamma(\mathbb{Z}_{2p}\times \mathbb{Z}_2)$ where $p\neq 2$ is not perfect.
    \item \label{sn-perfect} $\Gamma(S_n)$ is perfect if and only if $n \leq 3$. 
    \item     $\Gamma(A_n)$ is perfect if and only if $n \leq 5$. 
    \item \label{dihedral-perfect}    $\Gamma(D_n)$ is not perfect if and only if $n$ is divisible by at least two distinct odd primes.
\end{enumerate}  
\end{theorem}
\pf 
\begin{enumerate}
    \item Since $\mathbb{Z}_{pq}\cong \mathbb{Z}_p \times \mathbb{Z}_q$, the elements of order $p$ in $\mathbb{Z}_{pq}$ are $(1,0),(2,0),\ldots ,(p-1,0)$ and the elements of order $q$ are $(0,1),(0,2),\ldots,(0,q-1)$. Consider the following induced $5$-cycle $C$ in $\Gamma(\mathbb{Z}_{pq})$:  $$C:(p_1,0)\sim (0,0)\sim (-p_1,0)\sim (p_1,q_1)\sim (-p_1,q_2) \mbox{ where } p_1 \neq 0; q_1+q_2\neq 0. $$ Hence, $\Gamma(\mathbb{Z}_{pq})$ is not perfect.

    \item The elements of the group $\mathbb{Z}_{2p}\times \mathbb{Z}_2\cong \mathbb{Z}_2\times \mathbb{Z}_2\times \mathbb{Z}_p$ are of the form $(0,0,a),(0,1,a),(1,0,a),(1,1,a)$ where $a\in \mathbb{Z}_p$. Consider the following induced $5$-cycle $C$ in $\Gamma(\mathbb{Z}_{2p}\times \mathbb{Z}_2)$: $$C: (0,0,a)\sim (0,0,0)\sim (0,1,0)\sim (0,1,a)\sim (1,0,-a)\sim (0,0,a).$$ Hence $\Gamma(\mathbb{Z}_{2p}\times \mathbb{Z}_2)$ is not perfect.

    \item It can be checked that $\Gamma(S_2)$ and $\Gamma(S_3)$ are perfect. $\Gamma(S_4)$ contains the following induced $5$-cycle $C:e \sim (123)\sim (1342) \sim (1423) \sim (132)\sim e$ and hence not perfect. As $S_4$ is a subgroup of $S_n$ for $n\geq 5$, we get that $\Gamma(S_n)$ is not perfect for $n \geq 4$.

    \item It can be checked that $\Gamma(A_2),\Gamma(A_3),\Gamma(A_4)$ and $\Gamma(A_5)$ are perfect. As $S_n$ is a subgroup of $A_{n+2}$, we get that $\Gamma(A_n)$ is not perfect for $n \geq 6$.

    \item If $n$ is divisible by two distinct odd primes, say $p$ and $q$, then $D_n$ has a subgroup of isomorphic to $\mathbb{Z}_{pq}$ and hence $\Gamma(D_n)$ is not perfect. Conversely, if $n$ is not divisible by two distinct odd primes, then $n$ must be of the form $2^k$ or $2^kp^n$ or $p^n$, where $p$ is an odd prime. It is enough to show that $\Gamma(D_{2^kp^n})$ is perfect as $D_{2^k}$ and $D_{p^n}$ are subgroups of $D_{2^kp^n}$.
$$D_{2^kp^n}=\langle a,b: \circ(a)=2^kp^n, \circ(b)=2, ba=a^{-1}b\rangle.$$
We start by noting that \begin{itemize}
    \item any element of $D_{2^kp^n}$ is of the form $a^r$ or $a^rb$, where $0\leq r\leq 2^kp^n -1$.
    \item $a^r\sim a^sb$ for all $0\leq r,s\leq 2^kp^n -1$.
\end{itemize} 
Let $C:x_1\sim x_2\sim \cdots x_{2t+1}\sim x_1$ be an induced odd cycle in $\Gamma(D_{2^kp^n})$. If $x_1=a^rb$, then all the vertices, possibly except $x_2$ and $x_{2t+1}$, must be of the form $a^sb$. Now, if any one of $x_2$ or $x_{2t+1}$ is of the form $a^s$, then we get a chord in $C$, a contradiction. Thus all the vertices in $C$ are of the form $a^r$. Hence $C$ is an induced cycle in $\Gamma(\langle a\rangle)$, i.e., $\Gamma(\mathbb{Z}_{2^kp^n})$. But this can not happen as $\Gamma(\mathbb{Z}_{2^kp^n})$ is perfect (by Theorem \ref{abelian-perfect}). Hence $\Gamma(D_{2^kp^n})$ has no induced odd cycle of length $\geq 5$. Similarly, it can be shown that $\Gamma(D_{2^kp^n})^c$ has no induced odd cycle of length $\geq 5$. Thus $\Gamma(D_{2^kp^n})$ is perfect.
\end{enumerate}\qed 

\begin{theorem}
    If $p$ is an odd prime, then
    \begin{enumerate}
        \item $\Gamma(SL(2,p))$ and $\Gamma(GL(2,p))$ are not perfect.
        \item $\Gamma(PSL(2,p))$ is not perfect for $p\equiv \pm 1~(mod~8)$.
    \end{enumerate} 
\end{theorem}
\pf 
\begin{enumerate}
    \item The following is an induced $5$-cycle in $\Gamma(SL(2,p))$, and as $SL(2,p)\leq GL(2,p)$, both $\Gamma(SL(2,p))$ and $\Gamma(GL(2,p))$ are not perfect. The $p$'s over the edges denote the fact that order of the edge is $p$.
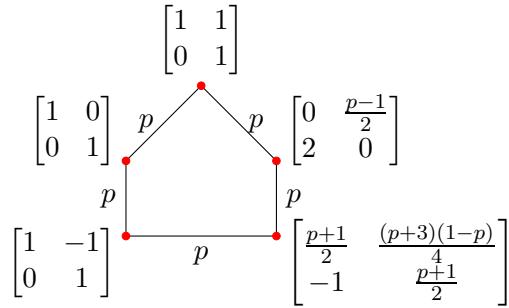
\begin{center}
    \begin{tikzpicture}[scale=0.5]
    {\small
   \draw [] (0, 0) -- (2, -2);
   \draw [] (0,0) -- (-2,-2);
   \draw [] (-2,-4) -- (-2,-2);
   \draw [] (2,-2) -- (2,-4);
   \draw [] (-2,-4) -- (2,-4);
\node [above] at (0, 0) {$\left[\begin{matrix}
    1 & 1\\
    0 & 1 \\
\end{matrix}\right]$};
			\draw [fill, red] (0, 0) circle[radius = 1mm]; 
   \node [right] at (2, -1.2) {$\left[\begin{matrix}
       0 & \frac{p-1}{2}\\
       2 & 0\\
   \end{matrix}\right]$};
			\draw [fill, red] (2, -2) circle[radius = 1mm];
   \node [left] at (-2, -1.2) {$\left[\begin{matrix}
       1 & 0\\
       0 & 1\\
   \end{matrix}\right]$};
			\draw [fill, red] (-2, -2) circle[radius = 1mm];
   \node [left] at (-2, -4.7) {$\left[\begin{matrix}
    1 & -1\\
    0 & 1 \\
\end{matrix}\right]$};
			\draw [fill, red] (-2, -4) circle[radius = 1mm]; 
   \node [right] at (2, -4.7) {$\left[\begin{matrix}
    \frac{p+1}{2} & \frac{(p+3)(1-p)}{4}\\
    -1 & \frac{p+1}{2} \\
\end{matrix}\right]$};
			\draw [fill, red] (2, -4) circle[radius = 1mm]; 
   \node [right] at (1,-1) {$\textcolor{black}{p}$};
   \node [left] at (-1,-1) {$\textcolor{black}{p}$};
   \node [left] at (-2,-3) {$\textcolor{black}{p}$};
   \node [right] at (2,-3) {$\textcolor{black}{p}$};
   \node [below] at (0,-4) {$\textcolor{black}{p}$};   }
\end{tikzpicture}
\end{center}

\item If $p\equiv \pm 1 ~(mod~8)$, then $PSL(2,p)$ contains a subgroup isomorphic to $S_4$ (Theorem 6.26, \cite{suzuki}). Thus, by Theorem \ref{bunch}(\ref{sn-perfect}), $\Gamma(PSL(2,p))$ is not perfect.\qed 
\end{enumerate}

\begin{remark}\label{psl-remark}
    Both of $\Gamma(PSL(2,3))$ and $\Gamma(PSL(2,5))$, i.e., $\Gamma(A_4)$ and $\Gamma(A_5)$ are perfect. It has been observed through numerical computations in SAGEMATH \cite{sage} that $\Gamma(PSL(2,p))$ is not perfect for $p\geq 7$. However, to resolve it we need to prove it for the cases when $p\equiv \pm 3~(mod~8)$. 
\end{remark}

Now we emphasize on perfectness of $\Gamma(G)$ in case of general groups. We start by recalling a graph, namely {\it subgroup sum graph}, which has been recently introduced in \cite{Cameron-Raveendra-Tamizh} . 
\begin{definition}
    Let $G$ be a finite abelian group, written additively, and $H$ be a subgroup of $G$. The subgroup sum graph $\mathcal{G}_{G,H}$ is the graph with vertex set $G$, in which two distinct vertices $x$ and $y$ are joined if $x+y \in H\setminus \{0\}$.
\end{definition}

It was proved in Corollary 3 of \cite{Cameron-Raveendra-Tamizh}, that $\mathcal{G}_{G,H}$ is perfect. Interestingly, if $G$ is a finite abelian $p$-group and $H$ is the subgroup generated by elements of order $p$ in $G$, then the prime order element graph of $G$, i.e., $\Gamma(G)$ coincides with $\mathcal{G}_{G,H}$ and hence $\Gamma(G)$ is perfect. However, this technique fails if $G$ is non-abelian $p$-group as the above argument does not work. Next, we provide only a sufficient condition for a finite $p$-group (not necessarily abelian) $G$ such that $\Gamma(G)$ is perfect.

\begin{theorem}\label{perfect-G}
    Let $p$ be a prime and $G$ be a finite $p$-group such that order of product of two elements of order $p$ is either $p$ or $1$, i.e., $\langle S \rangle =S \cup \{e\}$, where $S$ is the set of elements of order $p$ in $G$. Then $\Gamma(G)$ is perfect.
\end{theorem}
\pf If possible, let $\Gamma(G)$ be not perfect i.e. there exists an induced odd cycle such that $x_1\sim x_2 \sim x_3 \sim \cdots x_{2t+1} \sim x_1$ of length $2t+1\geq 5$. Then, $\circ(x_1x_2)=\circ({x_2}^{-1}{x_3}^{-1})=\circ(x_3x_4)=\circ({x_4}^{-1}{x_5}^{-1})=\cdots \circ({x_{2t}}^{-1}{x_{2t+1}}^{-1})=p$. This implies $\circ(x_1x_2{x_2}^{-1}{x_3}^{-1}x_3x_4)=\circ(x_1x_4)=p \mbox{ or } 1$. If $\circ(x_1x_4)=p$, then we have a chord between $x_1$ and $x_4$ which is a contradiction. Hence, $\circ(x_1x_4)=1$ i.e. $x_1x_4=e$ i.e. $x_4={x_1}^{-1}$. 

If $t=2$, then similarly, we get $x_5={x_2}^{-1}$, $x_3={x_1}^{-1}$ i.e. $x_3=x_4$ which is a contradiction. If $t\geq 3$, then we get $x_5={x_2}^{-1}$, $x_6={x_3}^{-1}$, $x_7={x_4}^{-1}$ i.e. $x_7=x_1$ which is again a contradiction. Hence $\Gamma(G)$ is perfect.\qed 

\begin{remark}
    Note that the converse of the above theorem is not true. For instance, we consider $G=(\mathbb{Z}_9\times \mathbb{Z}_3)\rtimes \mathbb{Z}_3$ [GAP Id:$(81,9)$], then $\Gamma(G)$ is perfect but $\langle S \rangle \neq S \cup \{e\}$.
\end{remark}

\begin{corollary}
    Let $G$ be a regular $p$-group. Then $\Gamma(G)$ is perfect.
\end{corollary}
\pf It follows from the fact that the condition $\langle S \rangle =S \cup \{e\}$ always hold in a regular $p$-group.

\begin{corollary}
    If $G$ is an abelian $p$-group, then $\Gamma(G)$ is perfect.
\end{corollary}
\pf It follows from the fact that abelian $p$-groups are regular.

Now, we prove a necessary and sufficient condition for an odd order group $G$ such that $\Gamma(G)$ is perfect.

\begin{theorem}\label{perfect_G-iff}
    Let $G$ be a group of odd order and $S$ denote the set of all prime order elements in $G$. Then $\Gamma(G)$ is perfect if and only if $\Gamma(\langle S\rangle)$ is perfect.
\end{theorem}
\pf As $\langle S\rangle$ is a subgroup of $G$, $\Gamma(G)$ is perfect implies $\Gamma(\langle S\rangle)$ is perfect. Conversely, let $\Gamma(\langle S\rangle)$ be perfect. If possible, let $C:x_1\sim x_2\sim \cdots \sim x_{2t+1}\sim x_1$ be an induced odd cycle of length $\geq 5$ in $\Gamma(G)$. Consequently, $$\circ(x_1x_2),\circ(x^{-1}_2x^{-1}_3),\circ(x_3x_4),\ldots,\circ(x_{2t+1}x_1) \mbox{ are all prime}.$$
Then $x^2_{1} \in \langle S\rangle$. Thus, $\circ(x_1\langle S\rangle)=2$ or $1$ in $G/\langle S\rangle$. However, as $G/\langle S\rangle$ is an odd order group, we must have $x_1\langle S\rangle=\langle S\rangle$, i.e., $x_1\in \langle S\rangle$. As vertices in $\langle S\rangle$ forms a component in $\Gamma(G)$, we must have $x_i \in \langle S\rangle$. But this implies that $C$ lies in $\Gamma(\langle S\rangle)$, a contradiction as $\Gamma(\langle S\rangle)$ is perfect. Similarly, it can be shown that $\Gamma(G)^c$ also has no induced odd cycle of length $\geq 5$. Eventually $\Gamma(G)$ is perfect.\qed

\begin{remark}
    It is to be noted that the above result does not hold true for groups of even order in general. For example, take $G\cong (\mathbb{Z}_2\times Q_8)\rtimes \mathbb{Z}_2$ (GAP id: $(32,44)$). Then $\langle S\rangle \cong (\mathbb{Z}_2\times \mathbb{Z}_4)\rtimes \mathbb{Z}_2$. Now, $\Gamma(\langle S\rangle)$ is perfect, but $\Gamma(G)$ is not.
\end{remark}

\begin{remark}\label{even-perfect-iff}
    Nonetheless, it is easy to see that the above theorem holds also in case of an even order group $G$ where $|G/\langle S \rangle|$ is odd.
\end{remark}

In the subsequent theorem, we characterize abelian groups $G$ for which $\Gamma(G)$ is perfect.
\begin{theorem}\label{abelian-perfect}
    Let $G$ be an abelian group. Then $\Gamma(G)$ is perfect if and only if $G$ has no subgroup $H$ of the form $\mathbb{Z}_{pq}$ or $\mathbb{Z}_{2p}\times \mathbb{Z}_2$, where $p,q$ are distinct odd primes. 
\end{theorem}
\pf Let $G$ be a group with a subgroup $H$ of the form $\mathbb{Z}_{pq}$ or $\mathbb{Z}_{2p}\times \mathbb{Z}_2$. By Proposition \ref{previous-results}(4), $\Gamma(H)$ is an induced subgraph of $\Gamma(G)$ and by Theorem \ref{bunch} (\ref{pq-not-perfect}) and (\ref{2p*2-not-perfect}), $\Gamma(H)$ is not perfect. Hence, $\Gamma(G)$ is not perfect.

Conversely, let $\Gamma(G)$ be not perfect. If possible, let $G$ does not contain a subgroup $H$ such that $H\cong \mathbb{Z}_{pq}$ or $H\cong \mathbb{Z}_{2p}\times \mathbb{Z}_2$ for some odd prime $p$ and $q$. As $G$ is abelian, by structure theorem of finite abelian groups, $G$ must be either a $p$-group or isomorphic to $\mathbb{Z}_{2^m}\times P$, where $P$ is an abelian $p$-group for an odd prime $p$. By Corollary \ref{2mpn-abelian-perfect}, $\Gamma(G)$ is perfect. This is a contradiction. Hence, $G$ contains a subgroup $H$ such that $H\cong \mathbb{Z}_{pq}$ or $H\cong \mathbb{Z}_{2p}\times \mathbb{Z}_2$.\qed

In the following consecutive theorems, we study the perfectness of $\Gamma(G)$ when $G$ is nilpotent.
\begin{theorem}
    Let $G$ be a finite nilpotent group in which one of the following conditions holds: 
    \begin{enumerate}
        \item $|G|$ has at least two distinct odd prime factors.
        \item $|G|=2^kp^n$ and the Sylow $2$-subgroup of $G$ is neither cyclic or nor quaternion group.
    \end{enumerate}
    Then $\Gamma(G)$ is not perfect.
\end{theorem}
\pf As $G$ is nilpotent, it is the direct product of its Sylow subgroups. 
\begin{enumerate}
    \item In this case, as $|G|$ has at least two distinct odd prime factors, say $p,q$, $G$ contains a subgroup isomorphic to $\mathbb{Z}_{pq}$. Now, as $\Gamma(\mathbb{Z}_{pq})$ is an induced subgraph of $\Gamma(G)$ and $\Gamma(\mathbb{Z}_{pq})$ is not perfect, our claim is established.
    \item In this case, $G$ has a subgroup isomorphic to $\mathbb{Z}_2\times \mathbb{Z}_2\times \mathbb{Z}_p$, and hence the result follows immediately.\qed 
\end{enumerate}

\begin{theorem}
    Let $G$ be a finite $p$-group for an odd prime $p$ such that order of product of two elements of order $p$ is either $p$ or $1$. Then both $\Gamma(\mathbb{Z}_{2^n}\times G)$ and $\Gamma(Q_{2^n}\times G)$ are perfect.
\end{theorem}
\pf If possible, let $\Gamma(\mathbb{Z}_{2^n}\times G)$ is not perfect. Suppose there exists an induced odd cycle $C: (x_1,y_1) \sim (x_2,y_2) \sim (x_3,y_3) \sim \cdots \sim (x_{2t+1},y_{2t+1}) \sim (x_1,y_1)$ of length $\geq 5$ in $\Gamma(\mathbb{Z}_{2^n}\times G)$. Two elements $(a,b)$ and $(c,d)$ in $\Gamma(\mathbb{Z}_{2^n}\times G)$ are adjacent if and only if either $\circ(ac)=1$ i.e. $ac=e_1$ and $\circ(bd)=p$ or $\circ(ac)=2$ and $\circ(bd)=p$ i.e. $bd=e_2$. In the first case, we call the corresponding edge, a $p$-edge, and similarly in the later case, we call it a $2$-edge.\\
{\it Claim:} There exist no two consecutive $2$-edges in $C$.\\
{\it Proof of claim:} Since the group $\mathbb{Z}_{2^n}$ contains unique element of order $2$ so does the group $\mathbb{Z}_{2^n}\times G$. If possible, let $(x_1,y_1),(x_2,y_2),(x_3,y_3)$ are $3$ elements in the cycle such that $\circ((x_1,y_1)(x_2,y_2))=\circ((x_2,y_2)(x_3,y_3))=2$ i.e. $\circ(x_1x_2)=\circ(x_3x_2)=2$ i.e. $x_1=x_3$. Also, $y_1y_2=y_2y_3=e_2$ i.e. $y_1=y_3$. Therefore, $(x_1,y_1)=(x_3,y_3)$ which is a contradiction.\\
{\it Claim:} $C$ can not be have only $p$-edges.\\
{\it Proof of Claim:} If all the edges are $p$-edges in $C$ then $\circ(x_ix_{i+1})=1$ and $\circ(y_iy_{i+1})=p$ for all $i=1, \ldots, 2t$. Then we get an induced odd cycle in $\Gamma(G)$, which is not possible by Theorem \ref{perfect-G}.\\
{\it Claim:} $C$ can have exactly one $2$-edge and rest are $p$-edges.\\
{\it Proof of claim:} Suppose there are more than one $2$-edges in $C$. Without loss of generality, let two nearest $2$-edges on $C$ be $$((x_1,y_1),(x_2,y_2)) \mbox{ and }((x_l,y_l),(x_{l+1},y_{l+1})),$$ where $l\geq 3$ and all edges of $C$ between them are $p$-edges. If $l$ is odd, say $2k+1$, then we have $\circ((x_1,y_1)(x_2,y_2))=\circ((x_{2k+1},y_{2k+1})(x_{2k+2},y_{2k+2}))=2$ and $\circ((x_i,y_i)(x_{i+1},y_{i+1}))=p$ for $i=2,3,\ldots,2k$. Therefore, we have 
\begin{itemize}
    \item $\circ(x_1x_2)=\circ(x_{2k+1}x_{2k+2})=2$.
    \item $x_2=x_4=x_6=\cdots=x_{2k}$.
    \item $x_3=x_5=\cdots=x_{2k+1}=x^{-1}_2$.
    \item $y_1=y^{-1}_2$ and $y_{2k+2}=y^{-1}_{2k+1}$.
    \item $\circ(y_2y_3)=\circ(y_3y_4)=\cdots =\circ(y_{2k}y_{2k+1})=p$.
\end{itemize} 
Using the last condition and the given hypothesis, we have $$\circ(y_1y^{-1}_{2k+1})=\circ((y_1y^{-1}_3)(y_3y_4)(y^{-1}_4y^{-1}_{5})\cdots (y^{-1}_{2k}y^{-1}_{2k+1}))=p \mbox{ or }1.$$
If $\circ(y_1y^{-1}_{2k+1})=p$, then using the fourth condition above, we have $\circ(y_1y_{2k+2})=p$. Also, as $\mathbb{Z}_{2^n}$ has a unique element of order $2$, we have $x_{2k+2}=x^{-1}_1$. Thus, we get $(x_1,y_1)\sim (x_{2k+2},y_{2k+2})$, a chord in $C$, a contradiction.

If $\circ(y_1y^{-1}_{2k+1})=1$, i.e., $y_1=y_{2k+1}$, then $(x_{2k+1},y_{2k+1})=(x^{-1}_2,y_1)$. Now, from the first two $p$-edges of $C$ we get $\circ(y_1y_4)=\circ((y_1y^{-1}_3)(y_3y_4))=p$ or $1$. If $\circ(y_1y_4)=p$, then we get a chord $(x_4,y_4)\sim (x_{2k+1},y_{2k+1})$ in $C$, a contradiction. So we have $y_4=y^{-1}_1$. But this implies $(x_2,y_2)=(x_4,y_4)$, a contradiction.

Thus $C$ can not have two $2$-edges at an odd distance $l$. Similarly it can be shown that $C$ can not have two $2$-edges at an even distance. Thus $C$ has exactly one $2$-edge.

Now, from the above claims, the cycle $C$ is of the form $(x_1,y_1) \sim (x_2,{y_1}^{-1}) \sim ({x_2}^{-1},y_3) \sim (x_2,y_4) \sim ({x_2}^{-1},y_5) \sim \cdots \sim ({x_2}^{-1},y_{2t+1}) \sim (x_1,y_1)$ where $\circ(x_1x_2)=2$. Clearly, $\circ(x_1{x_2}^{-1})=1$ i.e. $x_1{x_2}^{-1}=e_1$ i.e. $x_1=x_2$ i.e. $\circ(x_1)=4$. Then the cycle becomes $$(x_1,y_1) \sim (x_1,{y_1}^{-1}) \sim ({x_1}^{-1},y_3) \sim (x_1,y_4) \sim ({x_1}^{-1},y_5) \sim \cdots \sim ({x_1}^{-1},y_{2t+1}).$$

Now, we show that such a cycle can not exist. We do this by subdividing it into two cases, viz, when $C$ is a $5$-cycle or when $C$ is an odd cycle of length $\geq 7$.

Let $C$ be an induced odd cycle of length $\geq 7$. Then $\circ(y^{-1}_1y_3)=\circ(y^{-1}_3y^{-1}_4)=\circ(y_4y_5)=p$, i.e., $\circ(y^{-1}y_5)=p$ or $1$. If it is $p$, then we get a chord between the second and the fifth vertices. Thus $y_1=y_5$. Similarly, using the third, fourth and fifth edges of $C$, we get $\circ(y_3y_6)=1$, i.e., $y_6=y^{-1}_3$, and using fourth, fifth and sixth edges of $C$, we get $y_7=y^{-1}_4$. Thus we get a chord between the second and seventh vertices in $C$, a contradiction. Similarly, it can be shown that $C$ can not be of length $5$.

Thus $\Gamma(\mathbb{Z}_{2^n}\times G)$ has no induced odd cycle of length $\geq 5$. Proceeding similarly, it can be shown that $\Gamma(\mathbb{Z}_{2^n}\times G)^c$ has no induced odd cycle of length $\geq 5$. Eventually, $\Gamma(\mathbb{Z}_{2^n}\times G)$ is perfect.

The proof for perfectness of $\Gamma(Q_{2^n}\times G)$ follows analogously, as both $\mathbb{Z}_{2^n}$ and $Q_{2^n}$ have unique element of order $2$.\qed

\begin{corollary}\label{2mpn-abelian-perfect}
    If $G$ is an abelian $p$-group for an odd prime $p$, then both $\Gamma(\mathbb{Z}_{2^n}\times G)$ and $\Gamma(Q_{2^n}\times G)$ are perfect.
\end{corollary}

In the next few results, we classify the perfectness of $\Gamma(G)$ based on the order of $G$.
\begin{theorem}
    Let $G$ be a group of order $pq$ where $p,q$ are primes. Then $\Gamma(G)$ is perfect if and only if $G\not\cong \mathbb{Z}_{pq}$.
\end{theorem}
\pf If $p=q$, then $G$ is an abelian $p$ group and hence $\Gamma(G)$ is perfect. If $p=2$ and $q$ is an odd prime, then $G \cong D_q$ or $\mathbb{Z}_{2q}$ and hence $\Gamma(G)$ is perfect by Theorem \ref{bunch}(\ref{dihedral-perfect}) and Theorem \ref{abelian-perfect}. If $G$ is a non-cyclic group of order $pq$ with $p>q$, then all non-identity element of $G$ is of prime order. So, any vertex in $\Gamma(G)$ has degree either $pq-1$ or $pq-2$ and as a result, it can not have an induced odd cycle of length $\geq 5$. Thus, $\Gamma(G)$ is perfect.

\begin{theorem}
    Let $G$ be a group of order $2p^2$, where $p$ is an odd prime. Then $G$ is perfect if and only if $G \not\cong \mathbb{Z}_p \times D_p$.
\end{theorem}
\pf There are $5$ non-isomorphic groups of order $2p^2$, namely $\mathbb{Z}_{2p^2},\mathbb{Z}_p \times \mathbb{Z}_{2p},D_{p^2},(\mathbb{Z}_p \times \mathbb{Z}_{p})\rtimes \mathbb{Z}_2$ and $\mathbb{Z}_p \times D_p$. The first three of these are either abelian or dihedral groups, and hence by Theorems \ref{abelian-perfect} and \ref{bunch}(\ref{dihedral-perfect}), the corresponding graphs are perfect. In $(\mathbb{Z}_p \times \mathbb{Z}_{p})\rtimes \mathbb{Z}_2$, each non-identity element is of prime order, i.e., each vertex has degree either $2p^2-1$ or $2p^2-2$, and hence can not have an induced odd cycle of length $\geq 5$ and thereby making the corresponding graph perfect.

For $\mathbb{Z}_p \times D_p$, we exhibit an induced odd cycle in the graph. Any element of $\mathbb{Z}_p \times D_p$ is of the form $([i],a^jb)$, where $j=0,1,2,\ldots,p-1$. The following is an induced $5$-cycle in $\Gamma(\mathbb{Z}_p \times D_p)$: $([p-1],b)\sim([1],e)\sim ([0],e)\sim ([p-1],e)\sim ([1],ab)\sim ([p-1],b)$.

\begin{theorem}
    Let $G$ be a group of order $2pq$, where $p,q$ are distinct odd primes with $p<q$. Then $G$ is perfect if and only if $G \cong \mathbb{Z}_2 \times (\mathbb{Z}_q \rtimes \mathbb{Z}_p)$.
\end{theorem}
\pf If $p\nmid q-1$, then there are $4$ non-isomorphic groups of order $2pq$, namely $\mathbb{Z}_{2pq},D_{pq},\mathbb{Z}_p\times D_q$ and $\mathbb{Z}_q\times D_p$. All of those have elements of order $pq$ and hence the corresponding graphs are not perfect due to Theorem \ref{bunch}(\ref{pq-not-perfect}).

If $p\mid q-1$, then in addition to the above $4$ groups, we have two more groups of order $2pq$, namely $\mathbb{Z}_q \rtimes \mathbb{Z}_{2p}$ and $\mathbb{Z}_2 \times (\mathbb{Z}_q \rtimes \mathbb{Z}_p)$. In the first case, we have 
$\mathbb{Z}_q \rtimes \mathbb{Z}_{2p}\cong \langle x,y: \circ(x)=q,\circ(y)=2p,yxy^{-1}=x^k, \mbox{ where }\circ(k)=2p \mbox{ in }\mathbb{Z}^*_q\rangle,$ and $e \sim y^{p-1}\sim y \sim xy^{-1}\sim y^{1-p}\sim e$ is an induced $5$-cycle in $\Gamma(\mathbb{Z}_q \rtimes \mathbb{Z}_{2p})$.

For the group $G \cong \mathbb{Z}_2 \times (\mathbb{Z}_q \rtimes \mathbb{Z}_p)$, one can show that neither $\Gamma(G)$ nor its complement has any induced odd cycle, and hence $\Gamma(G)$ is perfect. 
%We omit the details of this proof.

\begin{theorem}
    Let $G$ be a $p$-group of order $\leq p^4$, where $p$ is a prime. Then $\Gamma(G)$ is perfect if and only if $G$ is not isomorphic to the group with GAP id $(81,7)$.
\end{theorem}
\pf If $p\geq 5$, then $G$ is a regular $p$-group and hence $\Gamma(G)$ is perfect. Similar argument holds for groups of order $3,3^2,3^3$. An exhaustive search using SAGEMATH for the groups of order $3^4$ and $2$-groups of order $\leq 16$, reveals that except the groups with GAP id $(81,7)$, all others have perfect graphs.\qed

\section{When is $\Gamma(G)$ a cograph}\label{cograph}
In this section, we search those groups $G$ for which $\Gamma(G)$ is a cograph. We begin with a lemma which rules out some of the cases when it is not a cograph.
\begin{lemma}\label{cograph-lemma}
    Let $G$ be a finite group. If $G$ contains 
    \begin{enumerate}
        \item an element of order $p^2$ for an odd prime $p$, or
        \item an element of order $pq$, where $p,q$ are distinct primes,
    \end{enumerate}
     then $\Gamma(G)$ is not a cograph.
\end{lemma}
\pf If $G$ contains an element of order $p^2$, then $G$ has a subgroup isomorphic to $\mathbb{Z}_{p^2}$ and $\Gamma(\mathbb{Z}_{p^2})$ has $[1]\sim [p-1]\sim [p+1]\sim [p^2-1]$ as an induced $P_4$.

If $G$ contains an element of order $2p$ with $p\geq 5$ being an odd prime, then $G$ has a subgroup isomorphic to $\mathbb{Z}_2\times \mathbb{Z}_p$ and $\Gamma(\mathbb{Z}_2\times \mathbb{Z}_p)$ has $([0],[0])\sim ([0],[2])\sim ([1],[p-2])\sim ([1],[4])$ as an induced $P_4$. For $p=3$, $\Gamma(\mathbb{Z}_6)$ has $[2]\sim [0]\sim [4]\sim [5]$ as an induced $P_4$. Finally, if $p,q$ are both odd primes, then as $\Gamma(\mathbb{Z}_{pq})$ is not perfect, it is also not a cograph.\qed 

\begin{corollary}
    Let $G$ be a finite group such that $\Gamma(G)$ is a cograph. Then $G$ is an EPPO group.
\end{corollary}
\pf It follows immediately from the above lemma and hence omitted. \qed

\begin{remark}
    The converse of the above corollary is not true, because EPPO groups allow elements of order $p^2$ for any prime $p$. However, if $p\geq 3$, then it will not be a cograph. So, the class of finite groups $G$ for which $\Gamma(G)$ is a cograph is a proper subclass of the class of EPPO groups.
\end{remark}

In view of the above remark, we first study the $p$-groups for which the graph is a cograph.

\begin{theorem}\label{pgroup-cograph-exponent-ood-prime}
    Let $G$ be a $p$-group for some odd prime $p$. Then $\Gamma(G)$ is a cograph if and only if exponent of $G=p$.
\end{theorem}
\pf If $\Gamma(G)$ is a cograph, then it follows from Lemma \ref{cograph-lemma}. Now, we assume that  $G$ be a $p$-group with exponent $p$. If possible, $\Gamma(G)$ is not a cograph, i.e., it contains an induced $P_4$, say $x_1\sim x_2 \sim x_3 \sim x_4$. This implies $x_1x_3=x_1x_4=e$, i.e., $x_3=x_4$, a contradiction. \qed

\begin{theorem}
    Let $G$ be a non-cyclic abelian $2$-group. Then $\Gamma(G)$ is a cograph if and only if exponent of $G \leq 4$.
\end{theorem}
\pf Let $G$ be a non-cyclic abelian $2$-group with $exp(G)\leq 4$. If possible, let $\Gamma(G)$ has a induced $P_4$, namely $x_1\sim x_2\sim x_3\sim x_4$. Thus, we have $$\circ(x_1x_2)=\circ(x_2x_3)=\circ(x_3x_4)=2 \mbox{ and}$$
$$\circ(x_1x_3),\circ(x_2x_4),\circ(x_1x_4)=1 \mbox{ or }4.$$

{\it Claim:} $\circ(x_1x_3)=4$.\\
{\it Proof of Claim:} Suppose $\circ(x_1x_3)=1$. Then $x_1=x^{-1}_3$, $\circ(x_1)=\circ(x_3)=4$ and $x^3_3=x_1$ (as otherwise $x_1=x_3$). Thus $x_1x_4=x^3_3x_4=x^2_3(x_3x_4)$ being a product of two elements of order $2$, must be of order $1$ or $2$. As $x_1\not\sim x_4$, we have  $\circ(x_1x_4)=1$, i.e., $x_1=x^{-1}_4$ which implies $x_3=x_4$, a contradiction. Hence the claim holds true.

Now, as $G$ is commutative, $\circ(x_1x_2\cdot x_2x_3)=1$ or $2$. If $x_1x^2_2x_3=e$, i.e., $x_2^2(x_1x_3)=e$. As $\circ(x_1x_3)=4$, we have $\circ(x_2^2)=4$, i.e., $\circ(x_2)=8$, a contradiction. So, we must have $(x_1x_2\cdot x_2x_3)^2=e$, i.e., $(x_1x_3)^2x^4_2=e$, i.e., $(x_1x_3)^2=e$ ( since $exp(G)\leq 4$). This implies $\circ(x_1x_3)=1$ or $2$, contradicting the above claim. So $\Gamma(G)$ is a cograph.

Conversely, let $\Gamma(G)$ be a cograph and if possible, $exp(G)\geq 8$. Then there exists an element $a$ of order $8$ in $G$. Again, as $G$ is non-cyclic and abelian, $G\neq \langle a \rangle$ and $G\setminus \langle a \rangle$ has at least one element $b$ of order $2$. Thus, $H=\langle a,b \rangle \cong \mathbb{Z}_8 \times \mathbb{Z}_2$ is a subgroup of $G$. However, as $\Gamma(\mathbb{Z}_8 \times \mathbb{Z}_2)$ has an induced $P_4$, we get a contradiction. Hence $exp(G)\leq 4$.\qed

It is specifically noted that the problem of classifying non-abelian $2$-groups $G$ for which $\Gamma(G)$ is a cograph, is still unresolved. We discuss a few open questions related to this in Section \ref{conclusion}. So, we close this section with the following theorem.

\begin{theorem}\label{cograph-final-theorem}
    Let $G$ be a finite group. If $\Gamma(G)$ is cograph, then the following holds:
    \begin{enumerate}
        \item The order of any nonidentity element of $G$ is either prime (possibly even) or a power of $2$.
        \item If $p$ is an odd prime dividing $|G|$, then Sylow $p$-subgroup of $G$ is of exponent $p$.
        \item If $H$ is a non-cyclic abelian Sylow $2$-subgroup of $G$, then exponent of $H\leq 4$.
        \item No two elements of distinct prime order commute.
        \item If $G$ is not a $p$-group, $Z(G)$ is trivial.
    \end{enumerate}
\end{theorem}
\pf The proof follows from the above two theorems and the properties of EPPO groups.\qed

It is to be noted that even if $G$ is a group satisfying all the five conditions, $\Gamma(G)$ may fail to be a cograph, e.g., $S_4$. A minor corollary follows from the above theorem.
 
\begin{corollary}
    Let $G$ be a group of order $m$ or $2m$, where $m$ is an odd integer, such that $\Gamma(G)$ is a cograph. Then every element of $G$ is of prime order.
\end{corollary}

\section{When is $\Gamma(G)$ chordal}\label{chordal-section}
In this section, we study those groups $G$ for which $\Gamma(G)$ is a chordal, interval, split and threshold graphs. We start with a lemma which will help us for characterizing the chordal graphs.
\begin{lemma}
    Let $G$ be a finite $2$-group such that $\Gamma(G)$ is chordal. Then $G$ is isomorphic to $\mathbb{Z}_{2^n}$ or $Q_{2^n}$ or $\mathbb{Z}^n_2$.
\end{lemma}
\pf If all elements of $G$ are of order $2$, then $G\cong \mathbb{Z}^n_2$ and $\Gamma(G)$ is complete and hence chordal. If $G$ has a unique element of order $2$, then $G$ is isomorphic to $\mathbb{Z}_{2^n}$ or $Q_{2^n}$ and $\Gamma(G)$ is a disjoint union of some copies of $K_2$ and isolated vertices, and thereby vacuously chordal.

So, we assume that $G$ has more than on element of order $2$ and at least one element of order $4$. Now, as $G$ is a $2$-group, $Z(G)$ is non-trivial.

{\it Case 1:} $Z(G)$ contains an element of order $4$, say $a$. Then $\circ(a^2)=2$. Thus there exist $b \in G \setminus \langle a \rangle$ such that $\circ(b)=2$. Then $H=\langle a,b\rangle\cong \mathbb{Z}_4\times \mathbb{Z}_2$.

{\it Case 2:} $Z(G)$ contains only elements of order $\leq 2$. In this case, we choose $a \in G \setminus Z(G)$ such that $\circ(a)=4$ and $b \in Z(G)\setminus \{a^2\}$ such that $\circ(b)=2$. Then $H=\langle a,b\rangle\cong \mathbb{Z}_4\times \mathbb{Z}_2$.

Now, as in both the cases, $G$ has a subgroup isomorphic to $\mathbb{Z}_4\times \mathbb{Z}_2$ and $\Gamma(\mathbb{Z}_4\times \mathbb{Z}_2)$ is not chordal, then we have $\Gamma(G)$ is not chordal, a contradiction. \qed 

\begin{theorem}
    $\Gamma(G)$ is chordal if and only if $G$ is isomorphic to one of the following groups: $S_3,\mathbb{Z}_3,\mathbb{Z}_{2^n},Q_{2^n},\mathbb{Z}^n_2$.
\end{theorem}
\pf Let $\Gamma(G)$ be chordal and $p\geq 5$ be an odd prime divisor of $|G|$ and $H$ be a subgroup of order $p$ in $G$. Let $x,y \in H$ such that $\circ(x)=\circ(y)=p$ and $x\neq y^{-1}$. Then $x \sim y\sim x^{-1}\sim y^{-1}\sim x$ is an induced $4$-cycle in $\Gamma(G)$ and hence $\Gamma(G)$ is not chordal. Thus, $|G|$ can not have a prime factor $\geq 5$.

If $9$ divides $|G|$, then $G$ has a subgroup $H$ isomorphic to $\mathbb{Z}_9$ or $\mathbb{Z}_3\times \mathbb{Z}_3$. In the first case, $\Gamma(H)$ has an induced $6$-cycle and in the later, it has an induced $4$-cycle. Thus, $\Gamma(H)$ and hence $\Gamma(G)$ can not be chordal, i.e., $9\nmid |G|$.

Hence $|G|=2^k,3\cdot 2^k$ or $3$. If $|G|=2^k$, by previous lemma, we get $G$ is isomorphic to $\mathbb{Z}_{2^n}$ or $Q_{2^n}$ or $\mathbb{Z}^n_2$. If $|G|=3$, then $G \cong \mathbb{Z}_3$. 

So we are left with the case when $|G|=3\cdot 2^k$. Note that $G$ does not have any element of order $6$, because $\Gamma(\mathbb{Z}_6)$ has an induced $4$-cycle and hence is not chordal. Thus all non-trivial elements of $G$ are of order $3$ or $2^r$, where $1\leq r \leq k$, i.e., $G$ is an EPPO group. Let $G_2$ and $G_3$ denote a Sylow $2$-subgroup and Sylow $3$-subgroup of $G$ respectively. Then $G_3\cong \mathbb{Z}_3$ and by previous lemma, $G_2$ is isomorphic to $\mathbb{Z}_{2^n}$ or $Q_{2^n}$ or $\mathbb{Z}^n_2$.

If $G_2\cong \mathbb{Z}_{2^n}$, then both Sylow subgroups of $G$ are cyclic and using Burnside's Transfer Theorem (Theorem 4.3 of Chapter 7, \cite{gorenstein}), it follows that $G_3\lhd G$. If $G_2\cong Q_{2^n}$, and then by Theorem 2.3(1) of \cite{eppo-shi-yang}, $G_3\lhd G$. Thus, in both the above cases, $G$ has a subgroup $HG_3$ of order $12$, where $H$ is a subgroup of $G_2$ of order $4$ (provided $n\geq 2$). Now, out of $5$ non-isomorphic groups of order $12$, none of them is chordal.  However, if $n=1$, i.e., $G\cong S_3$, then $\Gamma(G)$ is chordal.

If $G_2\cong \mathbb{Z}^n_2$, then $G$ is a solvable group with every non-trivial element of prime order. Then by main theorem of \cite{prime-order-elements-orginal}, $n\geq 2$ and $G_2\lhd G$, i..e, $G\cong \mathbb{Z}^n_2 \rtimes \mathbb{Z}_3$ and has a subgroup isomorphic to $\mathbb{Z}^2_2 \rtimes \mathbb{Z}_3\cong A_4$. As $\Gamma(A_4)$ is not chordal, $\Gamma(G)$ can not be chordal.

Thus, if $\Gamma(G)$ is chordal, then $G$ must be one of the five groups mentioned in the statement of the theorem. The converse is obvious.\qed 

\begin{corollary}
    $\Gamma(G)$ is chordal if and only if $\Gamma(G)$ is an interval graph.
\end{corollary}
\pf As interval graphs are always chordal, we assume that $\Gamma(G)$ is chordal. Then by the above theorem, $G$ is isomorphic to one of the groups: $S_3,\mathbb{Z}_3,\mathbb{Z}_{2^n}$, $Q_{2^n}$, $\mathbb{Z}^n_2$. The first two of them can be checked to be interval graphs. If $G\cong \mathbb{Z}_{2^n}$ or $Q_{2^n}$, then $\Gamma(G)$ is a disjoint union of $K_2$'s and some isolated vertices, and hence an interval graph. If $G\cong \mathbb{Z}^n_2$, then $\Gamma(G)$ is complete and hence an interval graph.\qed 

\begin{corollary}
    $\Gamma(G)$ is a split graph if and only if $G$ is isomorphic to one of the following groups: $\mathbb{Z}_3,\mathbb{Z}_4,S_3,Q_8,\mathbb{Z}^n_2$.
\end{corollary}
\pf Let $\Gamma(G)$ be a split graph. Then it is also a chordal graph and hence $G$ is isomorphic to one of the following groups: $S_3,\mathbb{Z}_3,\mathbb{Z}_{2^n},Q_{2^n},\mathbb{Z}^n_2$.

Now, as $\Gamma(\mathbb{Z}_{2^n})$ has two isolated vertices and rest is a disjoint union of $K_2$, from $\Gamma(\mathbb{Z}_{8})$ onwards, all of them are not split. By similar argument, it can be shown that among generalized quaternion groups, only $\Gamma(Q_{8})$ is split. \qed

\begin{corollary}
    $\Gamma(G)$ is a threshold graph if and only if it is a split graph.
\end{corollary}
\pf Since  a graph is threshold if and only if it is both cograph and split. The rsult follows from the above corollary.\qed

\section{When is $\Gamma(G)$ clawfree}\label{clawfree-section}
In this section, we study those groups $G$ for which $\Gamma(G)$ is clawfree. We begin with a result which provides a sufficient condition for a group $G$ such that $\Gamma(G)$ is clawfree. 
\begin{lemma}
    Let $G$ be a finite group. If $G$ contains 
    \begin{enumerate}
        \item an element of order $p^2$ for a prime $p\geq 5$, or
        \item an element of order $pq$, where $p,q$ are distinct primes,
    \end{enumerate}
     then $\Gamma(G)$ is not clawfree.
\end{lemma}
\pf If $G$ has an element $a$ of order $pq$ with $q>p$, then $e\sim a^p, e \sim a^{-p}, e \sim a^q$ forms a claw in $\Gamma(G)$. If $G$ has an element $a$ of order $p^2$ with $p\geq 5$, then $a\sim a^{p-1},a\sim a^{2p-1}$ and $a\sim a^{3p-1}$ forms a claw in $\Gamma(G)$. (It is to be noted that as $p\geq 5$, we have $3p-1<p^2$.)\qed

It is to be observed that from (2) of the above lemma, it follows that if $\Gamma(G)$ is  clawfree, then $G$ is an EPPO group. Also, if $p$ is an odd prime $\geq 5$ and $G$ is a $p$-group, then $\Gamma(G)$ is clawfree if and only if exponent of $G$ is $p$. The proof follows in a similar way as that of Theorem \ref{pgroup-cograph-exponent-ood-prime}. Thus, in the next two lemmas, we study the $2$-groups and $3$-groups only.

\begin{lemma}
    Let $G$ be a finite $3$-group. Then $\Gamma(G)$ is clawfree if and only if $G$ is a cyclic $3$-group or $exp(G)=3$.
\end{lemma}
\pf If $exp(G)=3$, i.e., any non-identity element of $G$ is of order $3$, then any two elements of $G$ are adjacent unless they are inverse to each other and hence $\Gamma(G)$ has no claw. Similarly, if $G$ is a cyclic $3$-group, it has exactly two elements of order $3$ and hence $\Gamma(G)$ is clawfree.

Conversely let $G$ be clawfree. If $G$ is cyclic, then we are done. On the other hand, suppose $G$ is a non-cyclic group. If possible, let $G$ has an element $a$ of order $9$. As $G$ is non-cyclic, $G$ has at least $8$ elements of order $3$. If $a\in Z(G)$, then we choose $b\in G\setminus \langle a \rangle$ such that $\circ(b)=3$. Then $a\sim a^2,a \sim a^5$ and $a\sim ba^2$ is a claw in $\Gamma(G)$. Suppose $a \notin Z(G)$. As $G$ is a  $3$-group, it has a non-trivial center. If $Z(G)$ contains an element $b$ of order $3$ (other than $a^3$ and $a^6$), then also $a\sim a^2,a \sim a^5$ and $a\sim ba^2$ forms a claw in $\Gamma(G)$. Finally we assume that $Z(G)$ contains exactly two elements of order $3$, namely $a^3$ and $a^6$, then also we get the same claw. Thus, $G$ can not have any element of order $9$, i.e., $exp(G)=3$.\qed

\begin{lemma}
    Let $G$ be a finite $2$-group. Then $\Gamma(G)$ is clawfree if and only if exactly one of the following holds:
    \begin{enumerate}
        \item $G$ is a cyclic $2$-group or a quaternion group.
        \item $G$ is non-cyclic, non-quaternion group and has no subgroup isomorphic to $\mathbb{Z}_8$ or $D_4$.
    \end{enumerate}
\end{lemma}
\pf If $G$ is a cylic $2$-group or a quaternion group, then it is clear that $\Gamma(G)$ is clawfree, as $G$ has a unique element of order $2$. So we assume that $G$ is non-cyclic, non-quaternion group and has no subgroup isomorphic to $\mathbb{Z}_8$ or $D_4$. Then $G$ has at least three elements of order $2$ and any non-identity element of $G$ is of order $2$ or $4$. Moreover, any two elements of order $2$ commute, as otherwise they will generate a dihedral subgroup of $G$, which is a contradiction. If possible, let $\Gamma(G)$ has a claw $x,y_1,y_2,y_3$ with $x\sim y_i$ and $y_i\not\sim y_j$. Thus $\circ(xy_i)=2$ and $\circ(y_iy_j)\neq 2$ for $i\neq j$.

{\it Claim 1:} At least two $y_i$'s must be order $4$.\\
{\it Proof of Claim 1:} If $\circ(y_i)=\circ(y_j)=2$, then $\circ(y_iy_j)=2$ and hence $y_i\sim y_j$, a contradiction. That means at most one $y_i$ can be of order $2$. If any $y_i=e$, then the order of other two $y_j$'s will be $4$ (by non-adjacency condition) and the claim holds. If there is no $y_i=e$, then also at least two of the $y_i$'s must be of order $2$. Hence Claim 1 holds. Without loss of generality, let $\circ(y_1)=\circ(y_2)=4$.

{\it Claim 2:} $\circ(x)=4$.\\
{\it Proof of Claim 2:} As $\circ(y_1)=4$ and $x\sim y_1$, we have $x\neq e$. If possible, let $\circ(x)=2$, then $y_1y_2=(y_1x)(xy_2)$ is a product of two elements of order $2$ and hence $\circ(y_1y_2)=1$ or $2$. However, as $y_1\not\sim y_2$, we get $y_1y_2=e$, i.e., $y_1=y^3_2$. Again as $e=y_1y_2=(y_1x)(xy_2)$, we have $y_1x=xy_2$, i.e., $y_1=xy_2x^{-1}$, i.e., $y^3_2=xy_2x^{-1}$. Thus, $\langle x,y_2:\circ(y_2)=4; \circ(x)=2; y^3_2=xy_2x^{-1}\rangle \cong D_4$, a contradiction. Consequently $\circ(x)=4$.

{\it Claim 3:} $\circ(y_3)=4$.\\
{\it Proof of Claim 3:} As $\circ(x)=4$, we have $\circ(y_3)=2$ or $4$. If possible, let $\circ(y_3)=2$. Also we have $\circ(xy_3)=2$. Accordingly, $\circ(x)=\circ(xy^2_3)=\circ(xy_3\cdot y_3)=2$ (as $xy_3\neq y_3$), a contradiction. Thus, $\circ(y_3)=4$.

Thus, $y_1y_2=(y_1x)(x^2)(xy_2)$ is a product of three elements of order $2$ and hence $\circ(y_1y_2)=1$ or $2$. However, as $y_1\not\sim y_2$, $\circ(y_1y_2)=1$, i.e., $y_1y_2=e$. Similarly, $y_2y_3=y_1y_3=e$. But, this implies $y_1=y_2=y_3=e$, a contradiction. Thus $\Gamma(G)$ is clawfree.

Conversely, let $\Gamma(G)$ be clawfree. If $G$ is cyclic or quaternion, then we are done. Moreover, $G$ can not have a subgroup isomorphic to $D_4$ as $\Gamma(D_4)$ has a claw, namely $b\sim e, b \sim a, b \sim a^3$ where $D_4=\langle a,b: \circ(a)=4; \circ(b)=2; ba=a^{-1}b\rangle$. Suppose $G$ is non-cyclic and non-quaternion and if possible, $G$ has a subgroup isomorphic to $\mathbb{Z}_8$, i.e., there exists $a\in G$ with $\circ(a)=8$. As $G$ is neither cyclic nor quaternion, $G$ has at least three elements of order $2$. If $a\in Z(G)$, then we choose $b\in G\setminus \langle a \rangle$ such that $\circ(b)=2$. Then $a\sim a^3,a \sim a^3b$ and $a\sim a^7b$ is a claw in $\Gamma(G)$. If $Z(G)$ contains an element $b$ (other than $a^4$), then also the same claw works. Finally, we assume that $Z(G)$ contains exactly one elements of order $a^4$, namely $a^4$, then also we get the same claw. Thus, $G$ can not have any element of order $8$, i.e, $G$ has no subgroup isomorphic to $\mathbb{Z}_8$.\qed

From the above three lemmas, we get the following necessary condition for $\Gamma(G)$ to be clawfree.

\begin{theorem}\label{clawfree-final-theorem}
    Let $G$ be a finite group. If $\Gamma(G)$ is clawfree, then the following holds:
    \begin{enumerate}
        \item If $p$ divides $|G|$ and $p\geq 5$, then Sylow $p$-subgroup of $G$ is of exponent $p$.
        \item If $3$ divides $|G|$, then Sylow $3$-subgroup of $G$ is either cyclic or of exponent $3$.
        \item If $2$ divides $|G|$, then Sylow $2$-subgroup of $G$ is either cyclic or quaternion or it has no subgroup isomorphic to $D_4$ or $\mathbb{Z}_8$.
        \item No two elements of coprime order commute.
        \item If $G$ is not a $p$-group, $Z(G)$ is trivial.
        \item $G$ is either a $p$-group or non-nilpotent.
    \end{enumerate}
\end{theorem}

It is to be noted that even if $G$ is a group satisfying all the above conditions, $\Gamma(G)$ may fail to be a cograph, e.g., $(\mathbb{Z}_4\times \mathbb{Z}_4)\rtimes \mathbb{Z}_3$. However, using the above theorem, one can prove that $A_5$ is uniquely recognizable from $\Gamma(A_5)$, because among groups of order $60$,  $\Gamma(A_5)$ is only claw-free graph. In fact, it is the only EPPO group among groups of order $60$. Moreover, one more interesting result can be made out of $A_5$.
\begin{theorem}
    Let $G$ be a non-abelian simple group such that $\Gamma(G)$ is clawfree. Then $G\cong A_5$.
\end{theorem}
\pf Since $\Gamma(G)$ is clawfree, $G$ is an $EPPO$ group. Suzuki \cite{eppo-suzuki} showed that there are exactly eight finite non-abelian simple $EPPO$ groups, namely $PSL(2,q)$ with $q=4,7,8,9,17$, $PSL(3,4)$, $Sz(8)$ and $Sz(32)$. Among those $PSL(2,7)$, $PSL(2,9),PSL(2,17)$, $PSL(3,4)$ has subgroups isomorphic to $D_4$ and $Sz(32)$ has elements of order $25$, and hence they are ruled out. For $PSL(2,8)$ and $Sz(8)$, it can be checked using SAGE \cite{sage} that their graphs have claws. So, we are left with $PSL(2,4)\cong A_5$ and it yields a clawfree graph.

\section{Conclusion and Open Issues}\label{conclusion}
In this article, although we have proved a few characterization theorems related to forbidden subgraphs in $\Gamma(G)$ like chordal, inerval, split, threshold graphs etc., a lot more is yet to be done in this direction. The table below (Table \ref{comp_table}) shows a few examples to establish that neither clawfree-ness nor being a cograph is implied by one other. Finally, we conclude the article with a few open issues which one may find interesting. 

\begin{table}[h]
    \centering
    \begin{tabular}{|c|c|c|c|c|}
        \hline
        Group & Perfect  & Clawfree  & Cograph & Chordal \\
        \hline \hline
        $S_3$ or $\mathbb{Z}^n_2$   & True  & True   & True & True  \\
        \hline
        $D_p$, odd prime $p$   & True  & True   & True & False  \\
        \hline
        $D_{2^n}$   & True    & False   & True & False  \\
        \hline
        $\mathbb{Z}_{3^n}$   & True  & True  & False & False \\
        \hline
        $D_{p^2}$, odd prime $p$  & True  & False  & False & False \\
        \hline
    \end{tabular}
    \caption{Comparison of few groups w.r.t forbidden subgraphs}
    \label{comp_table}
\end{table}

As mentioned earlier in Remark \ref{psl-remark}, we pose the following problem:

\begin{problem}
    If $p$ is a prime such that $p\equiv \pm 3~(mod~8)$, then $\Gamma(PSL(2,p))$ is not perfect.
\end{problem}

The converse of Theorem \ref{cograph-final-theorem} is not true. It may be nice to find a necessary and sufficient condition on a group $G$ to have $\Gamma(G)$ a cograph. By this, we mean to:

\begin{problem}
    Identify the groups $G$ for which $\Gamma(G)$ is a cograph.
\end{problem}
In this regard, we have another observation done in SAGEMATH, which one may try to prove:

\begin{problem}
    Let $G$ be a non-abelian $2$-group such that $\Gamma(G)$ is a cograph, then exponent of $G$ is either $4$ or $|G|/2$.
\end{problem}

Similar question can also be raised for clawfree graphs, i.e., 
\begin{problem}
    Identify the groups $G$ for which $\Gamma(G)$ is clawfree.
\end{problem}

\section*{Statements and Declarations}
The second author acknowledges the funding of DST-FIST Sanction no. $SR/FST/MS-I/2019/41$ and DST-SERB-MATRICS Sanction no. $MTR/2022/000020$, Govt. of India. 

\subsection*{Data Availability Statements}
Data sharing not applicable to this article as no datasets were generated or analysed during the current study.

\subsection*{Competing Interests} The authors have no competing interests to declare that are relevant to the content of this article.

 %\section*{References}


\begin{thebibliography}{9999}
%\bibitem{p^4-paper} B.N. Al-Hasanat and A. Almazaydeh, On classification of groups of order $p^4$, where $p$ is an odd prime, International Journal of Mathematics and Computer Science, 17(2022), no. 4, 1569-1593.
\bibitem{brachter-kaja-1} J. Brachter and E. Kaja, Classification of non-solvable groups whose power graph is a cograph, Journal of Group Theory, Vol. 26(4), 2023.
\bibitem{brachter-kaja-2} J. Brachter and E. Kaja, On groups with chordal power graph, including a classification in the case of finite simple groups, Journal of Algebraic Combinatorics, Vol. 58, 1095-1124, 2023.
 \bibitem{Cameron-Raveendra-Tamizh} P.J. Cameron, R.R. Prathap and T.T. Chelvam, Subgroup sum graphs of finite abelian groups, Graphs and Combinatorics, Vol 38, Article 114, 2022.
\bibitem{Cameron-survey} P.J. Cameron, Graphs defined on groups, Int. J. Group Theory, Volume 11, Issue 2, pp. 53-107, 2022.
\bibitem{pallabi-cograph} P.J. Cameron, P. Manna and R. Mehatari, On finite groups whose power graph is a cograph, Journal of Algebra, Vol. 591, 59-74, 2022.
\bibitem{eppo-cameron} P.J. Cameron and N.V. Maslova, Criterion of unrecognizability of a finite group by its Gruenberg-Kegel graph, Journal of Algebra, Vol. 607(A), 186-213, 2022.
\bibitem{strong-perfect} M. Chudnovsky, N. Robertson, P. Seymour, R. Thomas, The strong perfect graph theorem, Annals of Mathematics, Vol. 164 (1), 51–229, 2006.
 \bibitem{prime-order-elements-orginal} M. Deaconescu, Classification of Finite Groups with all Elements of Prime Order, Proceedings of the American Mathematical Society, Vol. 106, No. 3, pp. 625-629, 1989.
 \bibitem{gorenstein} D. Gorenstein, Finite Groups, AMS Chelsea Publishing, 1968.
 \bibitem{higman} G. Higman, Finite Groups in Which Every Element Has Prime Power Order, Journal of the London Mathematical Society, s1-32(3), 335-342, 1957.
\bibitem{pallabi-forbidden} P. Manna, P.J. Cameron and R. Mehatari, Forbidden Subgraphs of Power Graphs, Electronic Journal of Combinatorics, Vol. 28(3), P3.4, 2021.
\bibitem{tapa-1} T. Manna, A. Das and B. Bhattacharya, Prime Order Element Graph of a Group, {\it to appear in} Ricerche Di Matematica. \url{https://doi.org/10.1007/s11587-024-00906-0}
\bibitem{rotman-book} J.J. Rotman, An Introduction to the Theory of Finite Groups, 4th Edition, Graduate Text in Mathematics, Springer, 1995. 
\bibitem{eppo-shi-yang} W. Shi and W. Yang, On finite groups with elements of prime power orders, arxiv 2020, \url{https://arxiv.org/abs/2003.09445}.
 \bibitem{suzuki} M. Suzuki, Group Theory I, Springer, 1982.
\bibitem{eppo-suzuki} M. Suzuki, On a class of doubly transitive groups, Ann. Math. 75:1 (1962), 105-145.
\bibitem{sage} W. Stein and others: Sage Mathematics Software (Version 7.3), Release Date: 04.08.2016, {\tt http://www.sagemath.org}.
\bibitem{west-graph-book} D.B. West, Introduction to Graph Theory, Prentice Hall, 2001.
\end{thebibliography}
\end{document}